\documentclass[11pt,reqno]{amsart}

\usepackage{commandearticle_new}

\usepackage[francais]{layout}

\makeatletter
\def\subsection{\@startsection{subsection}{2}%
	\z@{.5\linespacing\@plus.7\linespacing}{.25\linespacing}%
	{\normalfont\bfseries}}
\def\subsubsection{\@startsection{subsubsection}{3}%
	\z@{.5\linespacing\@plus.7\linespacing}{.25\linespacing}%
	{\normalfont\itshape}}
\makeatother
\graphicspath{{Images/}}

\newenvironment{propon}[1][Proposition]{\begin{trivlist}
\item[\hskip \labelsep {\bfseries #1}]}{\end{trivlist}}

\everymath{\displaystyle}
\begin{document}
\title[Third order equivalent equation for relative velocity lattice Boltzmann schemes]{Third order equivalent equation
 for the relative velocity\\ lattice Boltzmann schemes with one conservation law}
 
 \author{François Dubois}
\address[François Dubois]{CNAM Paris, d\'epartement de math\'ematiques, Univ Paris Sud, Laboratoire de math\'ematiques, UMR 8628, Orsay, F-91405} 
\email{Francois.Dubois@math.u-psud.fr}

\author{Tony Février}
\address[Tony Février]{Univ Paris Sud, Laboratoire de math\'ematiques, UMR 8628, Orsay, F-91405, Orsay, F-91405}
\email{Tony.Février@math.u-psud.fr}

\author{Benjamin Graille}
\address[Benjamin Graille]{Univ Paris Sud, Laboratoire de math\'ematiques, UMR 8628, Orsay, F-91405, CNRS, Orsay, F-91405}
\email{Benjamin.Graille@math.u-psud.fr}

\maketitle

\begin{abstract}
{\bf Equation équivalente au troisième ordre pour les schémas de Boltzmann sur réseau à vitesse relative à une loi de conservation.}
Nous étudions la précision formelle des schémas de Boltzmann sur réseau à vitesse relative.
Ils diffèrent des schémas de d'Humières puisque l'ensemble des moments pilotant la relaxation dépend d'un champ de vitesse arbitraire fonction du temps et de l'espace.
 Nous traitons de l'asymptotique des schémas à vitesse relative pour une loi de conservation :
nous établissons l'équation équivalente du troisième ordre pour un nombre arbitraire de dimensions et vitesses.
\end{abstract}

\begin{abstract}
We study the formal precision of the relative velocity lattice Boltzmann schemes.
They differ from the d'Humières schemes by their relaxation phase: it occurs for a set of moments parametrized by a velocity field function of space and time.
We deal with the asymptotics of the relative velocity schemes for one conservation law: the third order equivalent equation is exposed for an
arbitrary number of dimensions and velocities. 
\end{abstract}

\section*{Version française abrégée}

Dans cette contribution, nous présentons l'équation équivalente au troisième ordre 
des schémas de Boltzmann sur réseau à vitesse relative introduits dans \cite{Fev:2014:0,Fev:2014:1} (inspirés par le schéma cascade \cite{Geier:2006:0}) lorsqu'il n'y a qu'une loi de conservation.
Ces schémas à vitesse relative diffèrent des schémas de d'Humières \cite{dHu:1992:0} par leur phase de relaxation : l'ensemble des moments est maintenant paramétré par un certain champ de vitesse fonction du temps et de l'espace.

En dimension $d\in\N^{\star}$ d'espace, nous considérons $\mL$ un réseau cartésien de pas $\dx$. 
Le pas de temps $\dt$ est déterminé par une mise à l'échelle acoustique après spécification d'une vitesse 
$\lambda=\dx/\dt$.
Nous notons $\ddqq$ un schéma à $q$ vitesses 
\smash{$\vj\in\R^d$, $0\leq j\leq q{-}1$}, 
auxquelles sont associées les distributions de particules $\fj$, $0\leq j\leq q{-1}$.

Le schéma consiste à faire évoluer les distributions en deux étapes : la phase de collision puis la phase de transport.
En notant $\smash{(\Pk[0],\ldots,\Pk[q-1])}$, l'ensemble des polynômes définissant les moments, la matrice de passage des distributions à ces moments est définie par (\ref{eq:MatMu}) où $\utilde(\vectx,t)\in\R^d$ est le paramètre champ de vitesse relatif pour $\vectx\in\R^d$, $t\in\R$.
Le vecteur des moments $\smash{\vectmu=(\mku[0],\ldots,\mku[q-1])}$ est alors défini par (\ref{eq:ftomu})
et la phase de relaxation est donnée par (\ref{eq:relaxationu}), où $\mkueq$ est l'équilibre du moment $\mku$ et $\smash{\vects=(\sk[0],\ldots,\sk[q-1])\in\R^q}$, le vecteur des paramètres de relaxation. La phase de transport est donnée par (\ref{eq:transport}), après retour aux distributions. Le choix $\utilde=\vectz$ permet de retrouver le schéma de d'Humières.


Nous travaillons avec une seule loi de conservation sur la densité $\rho$ définie par (\ref{eq::rho}) : ainsi $\sk[0]=0$ et $\sk\neq0$, $1\leq k\leq q-1$.
L'objectif est la détermination de l'asymptotique du troisième ordre d'un schéma $\ddqq$ général. Nous utilisons pour cela la méthode des équations équivalentes, un développement formel initialement utilisé pour l'étude de schémas aux différences finies \cite{Yan:1968:0,Lerat:1974:0}. Cette méthode a été adaptée aux schémas de d'Humières dans \cite{Dub:2008:0,Dub:2009:0} et utilisée pour les schémas à vitesse relative à deux lois de conservation \cite{Fev:2014:0}. Elle s'effectue grâce à un développement de Taylor de l'équation de transport (\ref{eq:transport}) à petit pas de temps $\dt$. Pour cela, les distributions de particules sont supposées être les restrictions sur le réseau $\mL$ de fonctions suffisamment régulières. Les paramètres de relaxation sont supposés constants et l'équilibre, provenant d'un vecteur de distributions $\vectfeq=(\fjeq[0], \ldots, \fjeq[q-1])$ indépendant de $\utilde$, vérifie (\ref{eq:mueq}).

Afin d'alléger les notations, un indice latin dans les signes de sommation décrit implicitement les entiers de 0 à $q{-}1$, tandis qu'un indice grec décrit les entiers de 1 à $d$. Nous noterons par ailleurs 
$\Delta=\dt$. Sous ces hypothèses, nous communiquons le résultat suivant.

\begin{propon}
 L'équation de conservation sur la densité s'écrit au troisième ordre grâce à l'équation (\ref{eq:masse3_0}).\\
 La quantité de mouvement à l'équilibre $\vectqeq$ est donnée par $\vectqeq= \sum_{j=0}^{q-1} \vj\fjeq$,
%
 la dérivée particulaire $\dtj$ par $\operatorname{d}_t^j=\Dt+\smash{\sum_{\alpha}}\vjc{\alpha}\Dx[\alpha]$, $0\leq j\leq q-1,$
  %
le défault de conservation $\thk$ par $\thk=\sum_{j}\Miju[kj]\operatorname{d}_t^j  \fjeq,$ $0\leq k\leq q-1,$
%
 et le paramètre de Hénon $\sig[k]$ \cite{Henon:1987:0} par $\sig[k]=1/\sk-1/2$, pour $\sk\neq0$, $0\leq k\leq q-1.$
\end{propon}
\vspace{0.4cm}
Ce résultat s'obtient en développant successivement l'équation de transport (\ref{eq:transport}) aux ordres de zéro à trois. La composition par la matrice des moments donne le développement équivalent sur les quantités macroscopiques. La clé pour monter en ordre est l'utilisation de lemmes de transition \cite{Fev:2014:1}: ces lemmes sont les développements en $\Delta$ à différents ordres des moments non conservés par la collision.

Dans l'équation de conservation (\ref{eq:masse3_0}), nous remarquons que les termes d'ordre inférieur ou égal à un sont indépendants de $\utilde$ grâce à l'hypothèse (\ref{eq:mueq}) sur l'équilibre. Nous rappelons que ce constat est aussi vrai dans le cas à deux lois de conservation \cite{Fev:2014:0}. Le paramètre $\utilde$ apparaît dans le terme d'ordre deux : il est susceptible d'influer sur le comportement dispersif du schéma. Nous notons enfin l'importance du défaut de conservation : cette quantité, dépendante de l'équilibre, pilote les termes d'ordre un (diffusion) et deux (dispersion).

\section*{Introduction}

The analysis of consistency for the lattice Boltzmann schemes uses some Chapman-Enskog type formal expansions. This method, based on a multiscale expansion at small Knudsen number \cite{ChaCow:1939:0}, has been a useful tool for these schemes \cite{dHu:1992:0,Chen:1992:0}. In this contribution, we choose an alternative method: the equivalent equations development proposed in \cite{Yan:1968:0,Lerat:1974:0} for the study of finite difference schemes and adapted to the lattice Boltzmann framework in \cite{Dub:2008:0}. The third order equivalent equations of a general $d$ dimensions and $q$ velocities---$\ddqq$---d'Humières scheme \cite{dHu:1992:0} for one and two conservation laws have been determined with this method \cite{Dub:2009:0}. It has also been used to build high order schemes \cite{Dub:2009:2,Dub:2011:0}, to enforce isotropy to an arbitrary order \cite{Augier:2011:0,Augier:2014:0} and to study some vectorial lattice Boltzmann schemes for hyperbolic systems \cite{Gra:2014:0}.

This method has already been extended to the relative velocity $\ddqq$ scheme \cite{Fev:2014:0,Geier:2006:0} for two conservation laws. We adapt the work to one conservation law: the third order equivalent equation is derived. In a first section, we recall the characteristics of the relative velocity $\ddqq$ scheme and exhibit the differences with the d'Humières scheme. The second section introduces the main result of the paper, the third order expansion on the conserved variable. The steps leading to this formal development are described in a third section. The influence of the relative velocity is studied: we particularly locate the terms involving it. We also check the consistency with the asymptotics obtained for the d'Humières scheme \cite{Dub:2009:0}.

\section{The relative velocity $\ddqq$ scheme}

In this section, we recall the derivation for one conservation law of the relative velocity lattice Boltzmann schemes introduced in \cite{Fev:2014:0} and inspired by the cascaded scheme \cite{Geier:2006:0}.
We consider $\lattice$, a cartesian lattice in $d$ dimensions with a typical mesh size $\dx$. The time step $\dt$ is determined by the acoustic scaling after the specification of the velocity scale $\lambda\in\R$ by the relation $\dt=\dx/\lambda$.
For the scheme denoted by \ddqq, we introduce $\smash{\vectV=(\vj[0],\ldots,\vj[q-1])}$ a set of $q$ velocities of $\R^d$. We assume that for each node $\vectx$ of $\lattice$, and  each $\vj$ in $\smash{\vectV}$, the point $\smash{\vectx+\vj \dt}$ is also a node of the lattice $\lattice$.
The aim of the $\ddqq$ scheme is to compute a particle distribution 
$\smash{\vectf = (\fj[0],\ldots,\fj[q-1])}$
on the lattice $\lattice$ at discrete values of time. 

The scheme splits into two phases for each time iteration: first, the relaxation that is non linear and local in space, and second, the linear transport for which an exact characteristic method is used. In the framework of the relative velocity lattice Boltzmann schemes, the relaxation is written into a basis of moments depending on a velocity field $\utilde(\vectx,t)$, a given function of space and time. In the following, we note the velocity field $\utilde$ because it appears only during the collision that is local in space.
The matrix of moments $\MatMu$, which is supposed to be invertible, is defined by
\begin{equation}\label{eq:MatMu}
 \Miju[kj] = \Pk(\vj-\utilde), \qquad 0\leq k,j\leq q{-}1,
\end{equation}
where $(\Pk[0],\ldots,\Pk[q-1])$ are some polynomials of $\R[\vars{X}{1},\ldots,\vars{X}{d}]$. We choose the first polynomials as follows: $\Pk[0]=1$, $\Pk=\vars{X}{k}$, $1\leq k\leq d$. 
The moments $\vectmu=(\mku[0],\ldots,\mku[q-1])$ are then given by
\begin{equation}\label{eq:ftomu}
\vectmu = \MatMu \; \vectf,
\end{equation}
so that $\mku$ is the \kieme moment in the frame moving at the velocity $\utilde$.
The relaxation phase is diagonal in the shifted moments basis 
\begin{equation}\label{eq:relaxationu}
\mkue = \mku{+}\sk (\mkueq{-}\mku), \qquad 0\leq k\leq q{-}1,
\end{equation}
where $\mkueq$ is the \kieme moment at equilibrium and $\sk$, the relaxation parameter associated with the \kieme moment for $0\leq k\leq q{-}1$. The vector of the distribution functions at the equilibrium $\vectfeq$ is chosen independent of the velocity field $\utilde$ so that we have
\begin{equation}\label{eq:mueq}
 \vectmequ = \MatMu \vectfeq = \MatMu \MatMinvz \vectmeqz.
\end{equation}
 We choose to work in the case of one conservation law on a scalar variable named here ``density'' defined by
\begin{equation}\label{eq::rho}
\rho=\sum_j \fj.
\end{equation}
Consequently, the associated relaxation parameter $\sk[0]$ is null.
The particle distributions are then computed by 
\begin{equation}\label{eq:mtof}
\vectfe = \MatMinvu\vectmue.
\end{equation}
The transport phase reads
\begin{equation}\label{eq:transport}
 \fj(\vectx,t+\dt) = \fje(\vectx-\vj\dt,t), \qquad 0\leq j\leq q{-}1.
\end{equation}
Let us note that the d'Humières scheme is enclosed in this framework by taking $\utilde=\vectz$. 

\section{Third order equivalent equation}

A Taylor expansion method \cite{Dub:2009:0} is performed to derive the third order equivalent equation of the relative velocity $\ddqq$ scheme for one conservation law, with any general equilibrium. We explicit a formal development of the scheme when $\dt$ and $\dx$ go to zero with $\lambda=\dx/\dt$ and the relaxation parameters $\vects$ fixed. We obtain in this way a set of partial differential equations that are consistent with the scheme at different orders.

 In the following, we do not adapt the Boltzmann scheme to a particular partial differential equation---at first or second order---but we put in evidence the various operators hidden inside the algorithm. The reasoning consists in a formal development of the transport phase (\ref{eq:transport}) at small $\Delta t$, assuming that all particle distributions are the restrictions on the discretized space of a  sufficiently regular distribution function. The Taylor formula can thus be used as much as wanted.
 
 The result we want to emphasize is given by the proposition \ref{th:ordre3_0}: this is the third order equivalent equation for the relative velocity $\ddqq$ scheme with one conservation law.
We particularly prove that the second order equations are the same as the ones of the $\ddqq $ d'Humières scheme \cite{Dub:2009:0}: this independence results from our choice for the equilibrium (\ref{eq:mueq}). The relation (\ref{eq:mueq}) means that the approximated physics do not depend on the velocity field $\utilde$ up to the second order. The dependence only appears at the third order.
 In the following, we note $\Delta=\Delta t$. Unless otherwise specified, a sum over a greek parameter goes from $1$ to $d$ and over a latin parameter from $0$ to $q-1$.

  \begin{proposition}\label{th:ordre3_0}
The third order conservation equation on the density reads:
 \begin{multline}\label{eq:masse3_0}
\Dt \rho+\sum_{\beta} \Dx[\beta]((\qeq)^{\beta})-\Delta\sum_{\beta}\!\sig[\beta]\Dx[\beta](\thz[\beta])=-\Delta^2 \sum_{\beta,j,l\geqslant1}\sig[\beta]\sig[l]~\vjc{\beta}~\Dx[\beta](\dtj\big(\Mijinvu \thk[l]\big)\\
+\frac{\Delta^2}{12}\sum_{\beta,\gamma,j}\vjc{\beta}\vjc{\gamma}\Dxy[\beta\gamma](\dtjf)
+\frac{\Delta^2}{6}\sum_{\beta}\Dxy[\beta t](\thz[\beta])+{\rm \mathcal{O}}(\Delta^3),
 \end{multline}
  where the momentum equilibrium $\vectqeq$ is given by $\vectqeq=\sum_{j=0}^{q-1} \vj\fjeq$, the particular derivative $\dtj$ by $\operatorname{d}_t^j=\Dt+\smash{\sum_{\alpha}}\vjc{\alpha}\Dx[\alpha]$, $0\leq j\leq q-1,$ the conservation default $\thk$ by $\thk=\smash{\sum_{j}}\Miju[kj]\operatorname{d}_t^j  \fjeq,$ $0\leq k\leq q-1,$ and the Hénon's parameter $\sig[k]$ {\rm \cite{Henon:1987:0}} by $\sig[k]=1/\sk-1/2$, for $\sk\neq0$, $0\leq k\leq q-1.$
\end{proposition}
\vspace{0.2cm}
We observe that each term of order $l$ in $\Delta$ contains only space derivatives of order $l+1$.

 \section{Principal ideas for the proof of the main result}

 The detailed calculations of the proof of the proposition \ref{th:ordre3_0} are available in \cite{Fev:2014:1}. To obtain the equivalent equation up to any order $p\in\N^{\star}$, a multidimensional Taylor formula relative to $\Delta$ is applied to the transport phase (\ref{eq:transport}). 
This expansion is then written in the moments basis after multiplying by $\Miju[kj]$
\begin{equation}\label{eq:dvgenu}
\sum_{\substack{0\leq j\leq q-1, \\0\leq l\leq p}} \frac{\Delta^l}{l!}\Miju\partial^l_t\fj=\sum_{\substack{0\leq j\leq q-1, \\|\frak{a}|\leq p}} \frac{(-\Delta)^{|\frak{a}|}}{\frak{a}!}(\vj)^{\frak{a}}\Miju\partial_{\frak{a}}\fje+{\rm \mathcal{O}}(\Delta^{p+1}), \quad0\leq k\leq q-1,
\end{equation}
where $\frak{a}=(a_1,\ldots,a_{d})\in\N^d,~ (\vj)^{\frak{a}}=\smash{\prod_{\alpha=1}^{d}} (\vjc{\alpha})^{a_{\alpha}},~ \partial_{\frak{a}}=\partial_{a_1}\cdots\partial_{a_{d}},~ \frak{a}!=\smash{\prod_{\alpha=1}^{d}}a_{\alpha}!$ and $|\frak{a}|=\sum_{\alpha=1}^{d} a_{\alpha}$.

This relation develops the different moments according to the choice of the moment index $k$. Taking $k=0$ yields to the conservation of the density, and $k>1$ to the expansions of the non conserved moments.

To obtain the third order equivalent equation on the density, we proceed order by order. The first step thus consists in deriving the zeroth order corresponding to $p=0$ in (\ref{eq:dvgenu}). This gives the following proposition: its proof can be found in \cite{Dub:2008:0}.

\begin{lemme}[The particles stay close to the equilibrium]\label{th:ordre0dh} 
The particle distributions can be \break expanded as
 \begin{equation*}\label{eq:or0eq}
\fj=\fjeq+{\rm \mathcal{O}}(\Delta),\  \ \fje=\fjeq+{\rm \mathcal{O}}(\Delta),\quad 0\leq j\leq q-1. 
\end{equation*}   
\end{lemme}
\vspace{0.4cm}
These identities are assumed to be formally derivable as much as wanted: this point is used to obtain the first order equivalent equation on $\rho$ taking $p=1$ and $k=0$ in (\ref{eq:dvgenu}). This equation is the truncation of (\ref{eq:masse3_0}) at the first order.

%

%

The key to gain an order is to develop relatively to $\Delta$ the non conserved moments. To get the second order equivalent equation on $\rho$, these moments must be expanded at the second order. We present these expansions in the third order case, the second one being obtained by truncation. 

\begin{lemme}[Third order transition lemma]\label{th:transition2}
The non conserved moments before and after the collision read,
\begin{gather}
\mku=~\mkueq-\Delta \Big(\frac{1}{2}+\sigma_k\Big)\xi_k(\utilde,\Delta,\vectsigma)+{\rm \mathcal{O}}(\Delta^3),\quad  1\leq k\leq q-1,\label{eq:mnc1}\\
\mkue=~\mkueq+\Delta\Big(\frac{1}{2}-\sigma_k\Big)\xi_k(\utilde,\Delta,\vectsigma)+{\rm \mathcal{O}}(\Delta^3),\quad 1\leq k\leq q-1,\label{eq:mnc2}
\end{gather}
 with 
\begin{equation*}
\xi_k(\utilde,\Delta,\vectsigma)=\thk{-}\Delta\sum_{j,l\geqslant1}\sigma_l~\Miju[kj]~\dtj\big(\Mijinvu \thk[l]\big).
\end{equation*}
\end{lemme}	
\vspace{0.4cm}
Two relations are combined to obtain these identities: the expansion (\ref{eq:dvgenu}) at the desired order and the relaxation phase (\ref{eq:relaxationu}). They give immediately the second order truncation of the relations (\ref{eq:mnc1},\ref{eq:mnc2}). This second order version is the keypoint to get the complete third order (\ref{eq:mnc1},\ref{eq:mnc2}). 

To get the second order equivalent equation on the density, we choose $p=2$ and $k=0$ in (\ref{eq:dvgenu}). The calculation necessitates the second order expansion of the post collision momentum that is given by (\ref{eq:mnc2}) taking $1\leq k\leq d$.

To get the complete third order equivalent equation on $\rho$, the same reasoning is used an order further. The required third order development of the post collision momentum is obtained using (\ref{eq:mnc2}) with $1\leq k\leq d$. The second order truncation of (\ref{eq:masse3_0}) is then useful to characterize this development up to the third order, that leads to the proposition \ref{th:ordre3_0}.

Let us now recover the equation on the density for the d'Humières scheme \cite{Dub:2009:0} corresponding to $\utilde=\vectz$. We focus on the second order terms, the first order being independent of $\utilde$. The quantity
$$-\frac{\Delta^2}{6}\sum_{\beta}\Dxy[\beta t](\thz[\beta]),$$
is already present in \cite{Dub:2009:0}.
The second term, independent of the relaxation parameters, reads
$$\sum_{\beta,\gamma,j}\vjc{\beta}\vjc{\gamma}\Dxy[\beta\gamma](\dtjf)
=\sum_{\beta,\gamma,l}\Dxy[\beta\gamma](\Lauz{\beta\gamma}{l}\thz[l]).$$
There is an only term depending on $\sigmav$ in (\ref{eq:masse3_0}). Splitting the time and space derivatives of the particular derivative $\dtj$, it reads
\begin{equation*}
\sum_{\beta,j,l\geqslant 1}\sig[\beta]\sig[l]~\vjc{\beta}~\Dx[\beta](\dtj\big(\Mijinvz \thz[l]\big)=\sum_{\beta,l\geqslant 1}\sig[\beta]\sig[l]~\Dxy[\beta t](\thz[\beta])
+\sum_{\beta,\gamma,l\geqslant 1}\sig[\beta]\sig[l]~\Dxy[\beta \gamma](\Lauz{\beta\gamma}{l}\thz[l]),
\end{equation*}
where $\Lauz{\beta\gamma}{l}$ is the momentum velocity tensor defined by $\Lauz{\beta\gamma}{l}=\smash{\sum_j}\vjc{\beta}\vjc{\gamma}\Mijinvz[jl]$, $1\leq\beta,\gamma\leq d$, $0\leq l\leq q-1$.
We have thus recovered all the terms exhibited in \cite{Dub:2009:0}.

\bibliographystyle{plain}
\bibliography{Bibliographie}

\end{document}